\documentstyle[amscd,amssymb,verbatim,12pt,diagrams]{amsart}

\diagramstyle[PostScript=dvips]

\newcommand{\irr}{\operatorname{irr}}

\newcommand{\Def}{\operatorname{Def}}
\newcommand{\can}{\operatorname{can}}
\newcommand{\reg}{\operatorname{reg}}
\newcommand{\Isom}{\operatorname{Isom}}
\newcommand{\Aut}{\operatorname{Aut}}
\newcommand{\CH}{\operatorname{CH}}
\newcommand{\bC}{{\bf C}}
\newcommand{\bp}{{\bf p}}
\newcommand{\del}{\partial}
\newcommand{\red}{\operatorname{red}}

\newcommand{\OO}{{\cal O}}

\newcommand{\DD}{{\cal D}}

\newcommand{\TT}{{\cal T}}

\newcommand{\G}{{\Bbb G}}

\newcommand{\lan}{\langle}
\newcommand{\ran}{\rangle}

\newcommand{\CC}{{\cal C}}

\newcommand{\Spec}{\operatorname{Spec}}

\renewcommand{\P}{{\Bbb P}}

\newcommand{\si}{\sigma}

\newcommand{\Pic}{\operatorname{Pic}}

\newcommand{\de}{\delta}
\newcommand{\eps}{\epsilon}

\numberwithin{equation}{section}

\newtheorem{thm}{Theorem}[section]
\newtheorem{prop}[thm]{Proposition}
\newtheorem{lem}[thm]{Lemma}
\newtheorem{cor}[thm]{Corollary}
\newenvironment{rem}{\vspace{3mm}\noindent
{\bf Remark.}}{\vspace{3mm}}

\newenvironment{ex}{\vspace{3mm}\noindent
{\bf Example.}}{\vspace{3mm}}

\newcommand{\Pf}{\noindent {\it Proof}}

\newcommand{\ov}{\overline}

\newcommand{\rk}{\operatorname{rk}}
\newcommand{\ra}{\rightarrow}

\newcommand{\HH}{{\cal H}}

\newcommand{\ZZ}{{\cal Z}}
\newcommand{\XX}{{\cal X}}
\newcommand{\LL}{{\cal L}}

\newcommand{\Res}{\operatorname{Res}}

\renewcommand{\a}{\alpha}
\renewcommand{\b}{\beta}
\newcommand{\om}{\omega}
\newcommand{\De}{\Delta}
\newcommand{\la}{\lambda}

\newcommand{\C}{{\Bbb C}}

\newcommand{\Q}{{\Bbb Q}}

\newcommand{\Ga}{\Gamma}

\newcommand{\wt}{\widetilde}

\newcommand{\sub}{\subset}
\newcommand{\ed}{\qed\vspace{3mm}}
\newcommand{\bm}{{\bf m}}
\newcommand{\eff}{\operatorname{eff}}
\newcommand{\MM}{{\cal M}}

\newcommand{\Mtwoeff}{\MM^{\frac{1}{2},\eff}_g}
\newcommand{\Meffpt}{\MM^{\frac{1}{r},\bm,\eff}_{g,n}}
\newcommand{\Meffptpt}{\MM^{\frac{1}{r},\bm,\eff}_{g,n,s}}
\newcommand{\Mspinpt}{\MM^{\frac{1}{r},\bm}_{g,n}}

\newcommand{\Mgn}{\MM_{g,n}}
\newcommand{\Mgns}{\MM_{g,n+s}}
\newcommand{\Mgnbar}{\ov{\MM}_{g,n}}
 
\newcommand{\Mto}{\MM_{2,1}}
\newcommand{\Mtobar}{\ov{\MM}_{2,1}}
\newcommand{\Mtospin}{\MM^{\frac{1}{r},2}_{2,1}}
\newcommand{\Mtospinbar}{\ov{\MM}^{\frac{1}{r},2}_{2,1}}
\newcommand{\Motbar}{\ov{\MM}_{1,3}}
\newcommand{\Motspinbar}{\ov{\MM}^{\frac{1}{r},(2,r-1,r-1)}_{1,3}}
\newcommand{\JJ}{{\cal J}}
\newcommand{\Mtospintwo}{\MM_{3,1}^{\frac{1}{2},2}}

\newcommand{\Mteffpt}{\MM^{\frac{1}{r},\bm,\eff}_{3,n}}
\newcommand{\Mtwoeffpt}{\MM^{\frac{1}{r},\bm,\eff}_{2,n}}

\title{Moduli spaces of curves with effective $r$-spin structures}
\author{A. Polishchuk}
\thanks{Supported in part by NSF grant}

\begin{document}
\begin{abstract}
We introduce the moduli stack of pointed curves equipped with effective
$r$-spin structures: these are effective divisors $D$ such that
$rD$ is a canonical divisor modified at marked points. 
We prove that this moduli space is smooth and describe its connected
components.
We also prove that it always contains a component that projects birationally
to the locus $S^0$ in the moduli space of $r$-spin curves consisting
of $r$-spin structures $L$ such that $h^0(L)\neq 0$. Finally, we study the 
relation between the locus $S^0$ and Witten's virtual top Chern class.
\end{abstract}
\maketitle

\bigskip

\section{Introduction}

Let us fix integers $g\ge 1$, $r\ge 2$, $n\ge 0$ 
and a vector $\bm=(m_1,\ldots,m_n)$ of non-negative integers such that
$$m_1+\ldots+m_n+rd=2g-2$$
for some integer $d\ge 0$.
Consider the moduli space $\Meffpt$ {\it of effective $r$-spin structures}
parametrizing collections
$(C,D,p_1,\ldots,p_n)$, 
where $C$ is a (connected) smooth complex projective
curve of genus $g$, $D\sub C$ is an effective divisor of degree $d$,
$p_1,\ldots,p_n$ are (distinct) marked points on $C$ such that
\begin{equation}\label{conddiv}
\OO_C(rD+m_1p_1+\ldots+m_np_n)\simeq\om_C
\end{equation}
(see section \ref{calcsec} for the precise definition of the moduli stack).

Our main result is the following theorem.

\begin{thm}\label{mainthm} 
(a) The stack $\Meffpt$ is smooth of dimension $2g-2+d+n$.

\noindent
(b) If $d\ge 0$ then there exists a point $(C,D,p_1,\ldots,p_n)$ in $\Meffpt$ such
that $h^0(D)=1$.
\end{thm}

The proof of this theorem will occupy sections \ref{calcsec} and 
\ref{degsec}. The idea of the
proof of part (a) is very simple. First, we
show that the dimension of $\Meffpt$ at every point is $\ge 2g-2+d+n$ representing
$\Meffpt$ as a degeneracy locus. Then we prove that the dimension
of the tangent space to $\Meffpt$ at every point is $\le 2g-2+d+n$.
Part (b) is proved using a degeneration argument: we find the required
structure on some nodal curve and then prove that it can be smoothened.

In the case $d=0$ the moduli space $\Meffpt$ are closely related to the
spaces $\HH(m_1,\ldots,m_n)$ studied in \cite{KZ}. Recall that the latter 
spaces
parametrize pairs $(C,\om)$, where $\om$ is a nonzero holomorphic
$1$-form on a curve $C$ such that zeroes of $\om$ have given multiplicities
$(m_1,\ldots,m_n)$. The main result of \cite{KZ} is the complete description of
connected components of $\HH(m_1,\ldots,m_n)$ (see Theorems 1 and 2 of 
{\it loc.~cit.}). Usually these spaces have
one or two connected components---the only exceptions are 
the spaces $\HH(2g-2)$ and $\HH(g-1,g-1)$ for $g\ge 4$ that have three
connected components.
Using the results of \cite{KZ} 
one can immediately describe connected components of
our moduli spaces $\Meffpt$ (for arbitrary $d\ge 0$)
because of the following theorem that will be proved in section
\ref{compsec}. 

\begin{thm}\label{conncompthm} 
Connected components of $\Meffpt$ are in one-to-one
correspondence with connected components of $\HH(r,\ldots,r,m_1,\ldots,m_n)$,
where $r$ is repeated $d$ times.
\end{thm}

For the detailed list of connected components of $\Meffpt$ see 
Corollaries \ref{conncor1} and \ref{conncor2}.

Our interest in the moduli spaces $\Meffpt$ is because of their relation 
to the moduli spaces of $r$-spin curves.
Recall that an {\it $r$-spin structure of type} $\bm$
on a smooth $n$-pointed curve $(C,p_1,\ldots,p_n)$ is a line bundle $L$ on $C$
together with an isomorphism $L^{\otimes r}\simeq \om_C(-\sum_i m_ip_i)$.
The moduli spaces $\Mspinpt$ of $r$-spin curves (i.e., curves with $r$-spin structures), 
or rather their compactifications,
have been used in \cite{JKV} to construct a cohomological field theory for every $r\ge 2$.
The most nontrivial feature of the moduli spaces $\Mspinpt$ used in this construction
is the existence of certain canonical class $c^{\frac{1}{r}}\in \CH^{-\chi}(\Mspinpt)_{\Q}$,
where $\chi=d-g+1$ (see \cite{JKV}, \cite{PV}, \cite{P}). This class is called
{\it Witten's virtual top Chern class} because on the open subset where the
spin structure $L$ has no global sections it coincides with 
$c_{-\chi}(R\pi_*\LL)=(-1)^{\chi}\cdot c_{-\chi}(R^1\pi_*\LL)$, where $\LL$ is the universal $r$-spin structure  
on the universal curve $\pi:\CC\ra\Mspinpt$. 
Thus, the difference $c^{\frac{1}{r}}-c_{-\chi}(R\pi_*\LL)$ is supported on the locus
$S^0\sub\Mspinpt$ where $h^0(L)\neq 0$. 
The locus $S^0$ is closely related to our moduli space $\Meffpt$. Namely, locally
\footnote{The problem with the
existence of such a global morphism is due to the presence of automorphisms of $r$-spin structures even when
the underlying curve has no automorphisms.} 
(or on the level of
the coarse moduli spaces) we have a morphism from $\Meffpt$ to $\Mspinpt$ with the image $S^0$. 
Our main theorem implies that the codimension of $S^0$ is equal to $-\chi$ (which is also the
degree of the class $c^{\frac{1}{r}}$). Hence,
the difference $c^{\frac{1}{r}}-c_{-\chi}(R\pi_*\LL)$ is a linear combination of the
classes of irreducible components of $S^0$ of maximal dimension. 
Finding 
coefficients in this linear combination seems to be crucial
for better understanding of Witten's virtual top Chern class. Some computations in this 
direction are given in section \ref{compsec}. They suggest
that at least in the case $-\chi\le 1$ the answer is quite simple.

Witten's virtual top Chern class can be extended to the natural 
compactifation of $\Mspinpt$ constructed by Jarvis \cite{Jarvis}.
The same is true about the class $c_{-\chi}(R\pi_*\LL)$.
Unfortunately, the naive attempt to extend the locus $S^0\sub\Mspinpt$ to 
the compactification by considering the locus where $h^0(L)\neq 0$
for stable curves leads to components of larger dimension. Still, we believe that there should be
a nice formula for the difference $c^{\frac{1}{r}}-c_{-\chi}(R\pi_*\LL)$ on the compactified
moduli space involving some natural analogue of $S^0$.

Let us point out some easy corollaries of Theorem \ref{mainthm}.

\begin{cor}\label{maindimcor} 
Let $S_r^i(g,\bm)\sub\Mspinpt$ be the locus of $r$-spin curves 
$(C,L,p_1,\ldots,p_n)$ such that $h^0(L)\ge i+1$. Then codimension of $S^r_i(g,\bm)$ in $\Mspinpt$ is
$\ge -\chi(L)+i=g-1-d+i$.
\end{cor}

In the case $r=2$ and $n=0$ the above corollary says that the locus
of smooth curves on which there exists a theta-characteristic $L$ with $h^0(L)\ge i+1$ has
codimension at least $i$. In fact, Teixidor i Bigas showed in \cite{T1}
that this codimension is at least $2i-1$. It is plausible that similar analysis is applicable to
more general loci $S_r^i(g,\bm)$.

\begin{cor} The locus $S_r^0(g,\bm)\setminus S_r^1(g,\bm)\sub\Mspinpt$ 
is smooth of dimension $2g-2+d+n$ and is non-empty if $d\ge 0$.
\end{cor}

\noindent
{\it Convention}. Throughout this paper we work over $\C$.

\noindent
{\it Acknowledgment}. I am grateful to A.~Beilinson for communicating the proof of
Lemma \ref{cantanglem} to me. Also, I'd like to thank C.~Faber, J.~Harris, T.~Kimura and 
A.~Vaintrob for valuable discussions.
Parts of this paper were written during the author's visits to Max-Planck-Institut
f\"ur Mathematik in Bonn and the Institut des Hautes \'Etudes Scientifiques.
I'd like to thank these institutions for hospitality and support. 

\section{Dimension calculations}\label{calcsec}

Assume that $g$, $r$, $n$, $\bm$ and $d$ are fixed as in the introduction.
First, let us give a precise definition of the moduli space
$\Meffpt$. 
Let $\CC\ra\Mgn$ be the universal curve over $\Mgn$, $\CC^{(d)}$ be its
$d$-th relative symmetric power, $\JJ^{2g-2}$ 
the relative Jacobian of degree $2g-2$ over
$\Mgn$ and let $\si_r^{\bm}:\CC^{(d)}\ra\JJ^{2g-2}$
be the morphism sending $(C,D,p_1,\ldots,p_n)$ to $(C,\OO_C(rD+\sum_i m_ip_i))$.
Here by the relative Jacobian $\JJ^m$ of degree $m$ we mean the stack over $\Mgn$ such that
for a scheme $S$ the category $\JJ^d(S)$ has objects $(\bC,\xi)$, where $\pi:\bC\ra S$ is
a family of smooth curves of genus $g$, $\xi$ is an element in the relative
Picard group $\Pic(\bC/S)=H^0(S,R^1\pi_*\OO^*_{\bC})$ that restricts to an element of degree $m$ in the
Picard group of every curve in this family (cf. \cite{G}). 
We define $\Meffpt\sub\CC^{(d)}$ so that the following square is cartesian:
\begin{equation}\label{maindiagram}
\begin{diagram}
\Meffpt &\rTo& \CC^{(d)} \nonumber\\
\dTo &&\dTo^{\si_r^{\bm}}\\
\Mgn &\rTo^{c}& \JJ^{2g-2}
\end{diagram}
\end{equation}
where the morphism $c$ sends a curve $C$ to $(C,\om_C)$. 
In particular, $\Meffpt$ is a closed substack in $\CC^{(d)}$.
Note that if one uses the universal Picard stack of degree $2g-2$
instead of $\JJ^{2g-2}$ the resulting fibered product will
be a $\G_m$-torsor over $\Meffpt$. This $\G_m$-torsor parametrizes choices of an isomorphism
(\ref{conddiv}).
To estimate the dimension of $\Meffpt$ we will identify this substack
with the degeneracy locus associated
with certain morphism of vector bundles on $\CC^{(d)}$. The idea is that
the condition (\ref{conddiv})
for an effective divisor $D$ of degree $d$ on a curve $C$ with marked points
$(p_1,\ldots,p_n)$ is equivalent to the condition that the natural linear
map
\begin{equation}\label{mainmap}
H^0(C,\om_C)\ra H^0(C,\om_C/\om_C(-rD-\sum_i m_ip_i))
\end{equation}
has nonzero kernel. These maps constitute a morphism of vector
bundles $\phi:V_1\ra V_2$ on $\CC^{(d)}$, where $V_1$ is the Hodge
bundle with the fiber $H^0(C,\om_C)$, $V_2$ is the bundle with the
fiber given by the target of (\ref{mainmap}).
Now it is clear that $\Meffpt$ as a set coincides with the locus
$Z_g\sub\CC^{(d)}$ where the rank of $\phi$ is $<g$.

\begin{lem}\label{dimlem} 
For every point $x\in\Meffpt$ one has $\dim_x\Meffpt\ge 2g-2+d+n$.
\end{lem}

\Pf . From the definition of $Z_g$ as a degeneracy locus
we immediately get that the codimension of $Z_g$ in $\CC^{(d)}$
is at most
$$\rk V_2-g+1=rd+\sum_i m_i-g+1=g-1$$
at every point. This gives the required estimate.
\ed

Next, we turn to the study of the tangent spaces to $\Meffpt$.
Let us recall what are the tangent spaces to the relevant moduli spaces.
For $(C,p)\in\Mgn$ we have 
$$T_{(C,p)}\Mgn=H^1(C,\TT_C(-p_1-\ldots-p_n)),$$
where $\TT_C$ is the tangent sheaf to $C$;
here and below we use the abbreviation $p=(p_1,\ldots,p_n)$.
For a line bundle $L$ of degree $d$ we have
$$T_{(C,p,L)}\JJ^d=H^1(C,A_{L,p}),$$
where $A_L$ is the Atiyah algebra of $L$, i.e., the sheaf of
differential operators $L\ra L$ of order $\le 1$,
$A_{L,p}\sub A_L$ is the subsheaf of operators with the
symbol vanishing at all the points $p_1,\ldots,p_n$ (see \cite{Diaz}). 
Finally, for an effective divisor $D$ of degree $d$ in $C$ we have
$$T_{(C,p;D)}\CC^{(d)}=H^1(C,\TT_C(-D-p_1-\ldots-p_n)).$$

\begin{lem}\label{rootdeflem} 
Let $L$ be a line bundle on $C$ such that 
$L^r\simeq\om_C(-\sum_i m_i p_i)$ and let $s_0\in H^0(C,\om_C\otimes L^{-r})$
be the section corresponding to the natural embedding $L^r\ra\om_C$. 
Then there is a canonical isomorphism of sheaves
$$A_{\om_C,p}\simeq A_{L,p}$$
such that an operator $\del:\om_C\ra\om_C$ corresponds to the operator
$$\del':L\ra L: s\mapsto \frac{\del(s^r s_0)}{rs^{r-1}s_0}.$$
\end{lem}

\Pf . First, we have to check that the map $\del\mapsto \del'$ above is well defined.
We know the symbol of an operator $\del\in A_{\om_C,p_1,\ldots,p_n}$ vanishes
at all points $p_1,\ldots,p_n$. Therefore, if $\pi_i$ is a local equation of
$p_i$ then for a every regular $1$-form $\a$ near $p_i$ we have
$$\del(\pi_i^{m_i}\a)=\pi_i^{m_i}\del(\a)+m_i\pi_i^{m_i-1}(\si(\del)\cdot\pi_i)\a,$$
where $\si(\del)$ is the symbol of $\del$. Since $\si(\del)$ vanishes at $p_i$
we obtain that $\del$ preserves the subsheaf $L(-m_ip_i)\sub L$.
This implies that the above morphism
$$A_{\om_C,p}\ra A_{L,p}$$
is well defined. Since it preserves the symbols and reduces to the multiplication
by $1/r$ on operators of order $0$, it is an isomorphism.
\ed

I learned the proof of the following lemma from A.~Beilinson.

\begin{lem}\label{cantanglem} The tangent space to the morphism
$c:\Mgn\ra\JJ^{2g-2}$ at a point $(C,p)$ is the map
$$H^1(\TT_C(-\sum_i p_i))\ra H^1(A_{\om_C,p})$$
induced by the map
$$\TT_C(-\sum_i p_i)\ra A_{\om_C,p}:v\mapsto L_v,$$
where $v$ is a vector field vanishing at the marked points,
$L_v:\om_C\ra\om_C$ is its action on $1$-forms by the Lie derivative: 
$$L_v(\eta)=d\lan\eta,v\ran,$$
where $\eta$ is a $1$-form,
$d:\OO_C\ra\om_C$ is the de Rham differential.
\end{lem}

\Pf . Let $R$ be a local artinian ring. For a curve with marked points 
$(C,p)$ let us denote by $\Aut_R(C,p)$ 
the sheaf of $R$-automorphisms of $(C,p)$ deforming the
trivial automorphism: for every open subset $U\sub C$ the group
$\Aut_R(C,p)(U)$ consists of automorphisms of 
$U\times\Spec(R)$ over $\Spec(R)$ preserving the marked points in $U$
and reducing to the identity modulo the maximal ideal of $R$. 
We can also consider similar sheaf $\Aut_R(C,p;L)$ adding
a line bundle $L$ on $C$ to our data. 
One can identify the category $\Def_R(C,p)$ (resp., 
$\Def_R(C,p;L)$)
of $R$-deformations of $(C,p)$ (resp., of 
$(C,p,L)$)
with torsors over $\Aut_R(C,p)$ 
(resp., $\Aut_R(C,p;L)$):
to a deformation one associates the torsor of its isomorphisms with
the trivial deformation. Now we have a canonical functor
$\Def_R(C,p)\ra\Def_R(C,p;\om_C)$ 
that associates to a deformation of pointed curves 
$\CC=(\bC,\bp)$
the deformation $\CC^{\can}=(\bC,\bp;\om_{\bC/\Spec(R)})$. 
In particular, if 
$\CC_0=(\bC_0,\bp^0)$ is a trivial deformation 
then we have a map 
$$\Isom_R(\CC,\CC_0)\ra\Isom_R(\CC^{\can},\CC_0^{\can})$$
between the corresponding torsors, compatible with the
natural homomorphism 
$$a_R:\Aut_R(C,p)\ra\Aut_R(C,p;\om_C).$$
It follows that the second torsor is the push-forward of the first
torsor with respect to $a_R$. Therefore, the tangent map to the
morphism $c:\Mgn\ra\JJ^{2g-2}$ is just a map induced on $H^1$
by the homomorphism $a_{R_0}$ for $R_0=\C[t]/t^2$. The sheaf of groups
$\Aut_{R_0}(C,p)$ can be identified with the sheaf
of vector fields on $C$ vanishing at the marked points.
Similarly, $\Aut_{R_0}(C,p,\om_C)$ can be identified with
the sheaf of order-$1$ differential operators $\om_C\ra\om_C$
with the symbol vanishing at the marked points. The homomorphism
$a_R$ is induced by the infinitesimal action of vector fields on
$\om_C$ that is by Lie derivatives. This is exactly the assertion
we wanted to prove.
\ed

\begin{prop}\label{cokerprop} 
For every point $(C,D,p)\in\Meffpt$ 
the cokernel of the linear map
$$T_{(C,D,p)}\CC^{(d)}\oplus T_{(C,p)}\Mgn
\stackrel{(d\si_r^{\bm},dc)}{\ra} 
T_{(C,p,\om_C)}\JJ^{2g-2}$$
has dimension $1$.
\end{prop}

We need one simple lemma for the proof.

\begin{lem}\label{divdifflem} 
Let $D\sub C$ be an effective divisor. We can view $D$ as a subscheme in $C$ and consider
the corresponding reduced subscheme $D^{\red}\sub C$ as a divisor on $C$.
Then the natural map
$$\OO_C(-D^{\red})/\OO_C(-D)\ra \om_C/\om_C(-D+D^{\red})$$
induced by the de Rham differential $d:\OO_C\ra\om_C$
is an isomorphism.
\end{lem}

\Pf . This follows immediately from the fact that the derivative map
$f\mapsto f'$ induces an isomorphism $t\C[t]/t^n\C[t]\ra \C[t]/t^{n-1}\C[t]$.
\ed 

\noindent
{\it Proof of Proposition \ref{cokerprop}}. 
Set $E=p_1+\ldots+p_n$.
We have to study the cokernel of the map
$$H^1(\TT_C(-D-E))\oplus H^1(\TT_C(-E))\stackrel{(\a,\b)}{\ra} 
H^1(A_{\om_C,p_1,\ldots,p_n})$$
where $\a$ and $\b$ are tangent maps to $\si_r^{\bm}$ and $c$ respectively.
By the definition, the map $\a$ is equal to the composition 
$$H^1(\TT_C(-D-E))\stackrel{\a'}{\ra} H^1(A_{L,p})\wt{\ra} 
H^1(A_{\om_C,p})$$
where $L=\OO_C(D)$, $\a'$ is the tangent map to the natural morphism
$\CC^{(d)}\ra\JJ^d$ while the second arrow is induced by the isomorphism of
Lemma \ref{rootdeflem}. Let $s_1\in L=\OO_C(D)$ be the natural section vanishing on $D$.
Then the map $\a'$ is induced by the $\OO_C$-linear
morphism of sheaves 
$$\TT_C(-D-E)\ra A_{L,p}:
v\mapsto (s\mapsto v(s/s_1) s_1,$$
where $v$ is a tangent vector on $C$ vanishing at $D+E$, $s/s_1$ is viewed as a rational function with poles at $D$. 
This morphism fits into the exact sequence
$$0\ra\TT_C(-D-E)\ra A_{L,p}\ra L\ra 0$$
where the morphism $A_{L,p}\ra L$ sends an operator $\del:L\ra L$
to $\del(s_1)\in L$. Therefore, we have an exact sequence 
\begin{equation}\label{H1exseq}
H^1(\TT_C(-D-E))\stackrel{\a'}{\ra} H^1(A_{L,p})\ra H^1(L)\ra 0.
\end{equation}
On the other hand, by Lemma \ref{cantanglem} the map $\b$ is induced by the sheaf morphism
$$\TT_C(-E)\ra A_{\om_C,p}:v\mapsto (\eta\mapsto d\lan\eta,v\ran)$$
Now the exact sequence (\ref{H1exseq}) shows that the cokernel we are interested in 
coincides with the cokernel of the map induced on $H^1$ by the following composition
of sheaf morphisms:
$$\TT_C(-E)\ra A_{\om_C,p}\wt{\ra} A_{L,p}\ra L$$
where the middle arrow is the isomorphism of Lemma \ref{rootdeflem}. Using the explicit
description of the relevant morphisms we see that the composed morphism
$\TT_C(-E)\ra L$ sends $v$ to 
$$\frac{d\lan s_1^rs_0, v\ran}{rs_1^{r-1}s_0},$$
where $s_0\in\OO_C(\sum_i m_ip_i)\simeq \om_C\otimes L^{-r}$ is the natural section.
This morphism fits into the following morphism of exact sequences
\begin{equation}
\begin{diagram}[size=5em]
0 \rTo^{} &\TT_C(-E) &\rTo^{s_1^r s_0}&\OO_C(-E)&
\rTo &\OO_C(-E)/\OO_C(-rD-E_{\bm}-E)&\rTo 0 \nonumber\\
&\dTo &&\dTo{d} &&\dTo & \\
0\rTo & L &\rTo^{rs_1^{r-1}s_0}&\om_C&\rTo &\om_C/\om_C(-(r-1)D-E_{\bm})&\rTo 0
\end{diagram}
\end{equation}
where $E_{\bm}=m_1p_1+\ldots+m_np_n$. 
Hence, we get a morphism between the exact sequences of cohomology groups
\begin{equation}
\begin{diagram}
H^0(\OO_C(-E)/\OO_C(-rD-E_{\bm}-E))&\rTo &
H^1(\TT_C(-E)) &\rTo & H^1(\OO_C(-E))&\rTo 0 \nonumber\\
\dTo &&\dTo &&\dTo \\
H^0(\om_C/\om_C(-(r-1)D-E_{\bm})) &\rTo &H^1(L) &\rTo & H^1(\om_C)&\rTo 0
\end{diagram}
\end{equation}
Note that the map $H^1(\OO_C(-E))\ra H^1(\om_C)$ induced by the de Rham
differential is zero. On the other hand,
applying Lemma \ref{divdifflem} we find that the vertical arrow on the left is surjective.
Hence, the cokernel of the middle vertical arrow can be identified with the
$1$-dimensional space $H^1(\om_C)$.
\ed

\noindent
{\it Proof of part (a) of Theorem \ref{mainthm}.} By Lemma \ref{dimlem}
for every point $x\in\Meffpt$ one has 
$$\dim_x\Meffpt\ge 2g-2+d+n.$$
On the other hand, Proposition \ref{cokerprop} implies that
$$\dim T_x\Meffpt\le 2g-2+d+n.$$
Therefore, $\dim T_x\Meffpt=\dim_x\Meffpt=2g-2+d+n$, so $\Meffpt$
is smooth of this dimension.
\ed

\section{Degeneration argument}
\label{degsec}

\subsection{Construction of effective $r$-spin structures on
nodal curves}

Given the data $(g, r, n, \bm, d)$ as in the introduction with $d\ge 0$
we want to construct a nodal curve $C$ of arithmetic genus $g$
with $n$ smooth points $p_1,\ldots,p_n\in C$ and an effective divisor
$D$ of degree $d$ on the smooth part of $C$ such that $h^0(D)=1$ and
the divisor $rD+\sum_i m_ip_i$ is in the canonical linear series of $C$.
The first natural attempt would be to take $C$ rational (with $g$ nodes).
Let us formulate the corresponding problem about polynomials.

\vspace{2mm}

\noindent {\bf Problem}. For the data $(g,r,n,\bm,d)$ such that $d\ge 0$
find distinct complex numbers $x_1,\ldots,x_g;y_1,\ldots,y_g;z_1,\ldots,z_n$
and a polynomial $P\in\C[t]$ of degree $d$ such that $P(x_j)\neq 0$,
$P(y_j)\neq 0$ for all $j=1,\ldots,g$ and the following conditions hold:

\noindent
(i)$\Res_{t=x_j} R(t)dt+\Res_{t=y_j} R(t)dt=0$ for all $j=1,\ldots,g$, where
$$R(t)=\frac{P(t)^r\prod_{i=1}^n(t-z_i)^{m_i}}{\prod_{j=1}^g(t-x_j)(t-y_j)};$$ 

\noindent
(ii)if $f$ is a polynomial of degree $<d$ such that
$$\frac{f(x_j)}{P(x_j)}=\frac{f(y_j)}{P(y_j)}$$
for all $j=1,\ldots,g$ then $f=0$.

Indeed, the solution to this problem would give a curve $C$ with $g$ nodes obtained from $\P^1$
by identifying each $x_j$ with $y_j$. Furthermore, the points $z_i$ should be considered
as marked points. Then we claim that the divisor $D$ of zeroes of $P$ will have the
required properties. To see this we observe that condition (i) garantees that
the rational $1$-form $R(t)dt$ on $\P^1$ descends to the regular $1$-form on $C$
with the divisor of zeroes equal to $rD+\sum_i n_iz_i$. On the other hand, condition (ii)
implies that $H^0(C,\OO_C(D))=0$.

We don't know how to solve this problem in general. So our strategy will be first to
find the solution in the case $d\le 1$ and then to construct the nodal effective $r$-spin curve
$(C,D,p_1,\ldots,p_n)$ using these solutions (for $d>1$ the curve $C$ will be reducible).
Note that for $d\le 1$ and $g\ge 1$ the condition (ii) is satisfied automatically.
Moreover, the problem for $d=1$ reduces to a similar problem
for $d=0$ and with one more point added (the corresponding weight is $r$).
The following theorem gives a solution of the above problem for $d=0$.

\begin{thm}\label{P1thm} 
For every collection of non-negative integers $(m_1,\ldots,m_n)$ and
$\sum_i m_i=2g-2$ there exists distinct complex numbers $x_1,\ldots,x_g;y_1,\ldots,y_g;z_1,\ldots,z_n$ 
such that the rational function 
\begin{equation}\label{ratfuneq}
R(t)=\frac{\prod_{i=1}^n(t-z_i)^{m_i}}{\prod_{j=1}^g(t-x_j)(t-y_j)}
\end{equation}
satisfies $\Res_{t=x_j} R(t)dt+\Res_{t=y_j} R(t)dt=0$ for all $j=1,\ldots,g$.
\end{thm}

\Pf . The idea is to start with a degenerate solution of the equations of condition (i) and
then to deform them. The degenerate solution is obtained by taking $y_j=-x_j$ for
all $j$ and $z_i=0$ for all $i$. Let us fix distinct nonzero numbers $x_1,\ldots,x_g$ and
consider the variety $Z=Z(x_1,\ldots,x_g)$ consisting of the $(n+g)$-tuples $(z_1,\ldots,z_n,
y_1,\ldots,y_g)$ such that 
the corresponding function $R$ given by (\ref{ratfuneq}) satisfies (i) and in addition all points
$(y_1,\ldots,y_g)$ are distinct and disjoint from
the set $\{x_1,\ldots,x_g,z_1,\ldots,z_n\}$ (we do not require $z_i$'s to be distinct). 
We are going to study $Z$ near the point
$p_0=(0,\ldots,0,-x_1,\ldots,-x_g)$. Namely, we claim that for an appropriate choice of $x_1,\ldots,x_g$
the natural projection 
$$\pi:Z\ra\C^n:(z_1,\ldots,z_n,y_1,\ldots,y_g)\mapsto (z_1,\ldots,z_n)$$ 
is \'etale at $p_0$ (in particular, $Z$ is nonsingular at $p_0$).
This claim immediately implies the theorem. Indeed, we can take an $n$-tuple
of distinct numbers $z=(z_1,\ldots,z_n)$
in a sufficiently small neighborhood of $0\in\C^n$ and then find a point $p\in Z$ such that $\pi(p)=z$.
To prove the claim let us write explicitly the equations defining $Z$:
\begin{equation}\label{Zeq}
\frac{\prod_{i=1}^n(x_j-z_i)^{m_i}}{\prod_{k:k\neq j}(x_j-x_k)(x_j-y_k)}=
\frac{\prod_{i=1}^n(y_j-z_i)^{m_i}}{\prod_{k:k\neq j}(y_j-x_k)(y_j-y_k)}, \text{ for } j=1,\ldots,g.
\end{equation}
After some simplifications the differentials of these equations at $p_0$ can be written as follows:
\begin{equation}\label{tangenteq}
x_j\sum_{k:k\neq j}\frac{dy_j-dy_k}{x_j-x_k}+2\sum_{i=1}^n m_idz_i=(2g-2)dy_j, \text{ for } j=1,\ldots,g.
\end{equation}
To show that the projection $\pi$ is \'etale at $p_0$ it suffices to prove that after substituting $dz_i=0$
for all $i$ in (\ref{tangenteq}) we obtain a linear system for $dy_i$'s with only zero solution.
So we are reducing to showing that there exists $x_1,\ldots,x_g$ such that the linear system
on variables $t_1,\ldots,t_g$ 
\begin{equation}\label{lineareq}
x_j\sum_{k:k\neq j}\frac{t_j-t_k}{x_j-x_k}=(2g-2)t_j, \text{ for } j=1,\ldots,g
\end{equation}
has only zero solution. Let us prove this by induction in $g$. For $g=1$ the assertion is clear. Now let $g>1$.
By induction assumption we can
choose $x_2,\ldots,x_g$ (nonzero and distinct)
such that the above system for $g-1$ variables $t_2,\ldots,t_g$ has only zero solution.
Now we observe that if we fix these $x_2,\ldots,x_g$ and let $x_1$ tend to infinity then the system
(\ref{lineareq}) will tend to the following system:
\begin{align*}
&\sum_{k:k\neq 1}(t_1-t_k)=(2g-2)t_1,\\
&x_j\sum_{k:k\neq j,k\neq 1}\frac{t_j-t_k}{x_j-x_k}=(2g-2)t_j, \text{ for } j=2,\ldots,g
\end{align*}
By our choice of $x_2,\ldots,x_g$ the solution of such system necessarily has $t_2=\ldots=t_g=0$. Then
the first equation shows that $t_1=0$. 
\ed

Now we can construct the nodal $r$-spin curve with the required properties.

\begin{thm}\label{nodalthm} 
For every data $(g, r, n, \bm, d)$ as in the introduction such that $d\ge 0$, there exists
a rational nodal curve $C$ of arithmetic genus $g$, $n$ smooth distinct
points $p_1,\ldots,p_n\in C$ and an effective divisor
$D$ of degree $d$ supported on the smooth part of $C$, such that 
$\OO_C(rD+\sum_i m_ip_i)\simeq\om_C$ and $h^0(D)=1$.
\end{thm}

\Pf . Assume first that $d$ is even and $m_i=0$ for all $i$. 
Then our curve $C$ will be obtained by certain identifications from the disjoint union
$\sqcup_{s=1}^{d/2} (C_s\sqcup C'_s)$ of $d$ copies of $\P^1$ that are grouped in pairs $(C_s,C'_s)$.
Let us fix a pair of complex numbers $a\neq b$ such that $a^r=b^r$ and $|a|=|b|\neq 1$. Also, set
$\zeta_j=\exp(2\pi i j/r)$ for $j=1,\ldots,r$. We denote by $a(s),b(s),\zeta_j(s)$ (resp., $a'(s),b'(s),\zeta'_j(s)$)
these numbers considered as points on $C_s$ (resp., $C'_s$). Here are the identifications we need to make to get $C$:

\noindent (a) $a(s)$ is glued to $b(s)$ for all $s$;

\noindent (b) for $j=2,\ldots,r$; $s=1,\ldots,d/2$, each $\zeta_j(s)$ is glued to $\zeta'_j(s)$;

\noindent (c) for $s=1,\ldots,d/2-1$, $\zeta_1(s)$ is glued to $\zeta'_1(s+1)$; $\zeta_1(d/2)$ is glued to
$\zeta'_1(1)$.

Thus, $C$ is a kind of a wheel formed by $d/2$ curves $C_s\cup C'_s$, where $C_s$ and $C'_s$ each have one node
and are glued with each other in $r-1$ points.
We define the divisor $D$ by setting $D=\sum_s 0(s)+0'(s)$, where $0(s)$ (resp., $0'(s)$) is $0\in\P^1$ considered
as a point of $C_s$ (resp., $C'_s$). The presence of one node on each component garantees that $h^0(D)=1$.
To show that $\OO_C(rD)\simeq\om_C$ we observe that the rational $1$-form 
$$R(t)dt:=\frac{t^rdt}{(t^r-1)(t-a)(t-b)}$$
has opposite residues at $a$ and $b$ (since $a^r=b^r$). 
Hence, we can define a regular $1$-form $\eta$ on $C$ with the divisor $rD$
by setting $\eta|_{C_s}=R(t)dt$, $\eta|_{C'_s}=-R(t)dt$. 

Next, consider the case when $d$ is even but $m=\sum m_i>0$.
Then $C$ will be glued from the $d+1$ copies of $\P^1$ that are grouped as follows: $d/2$ pairs $(C_s,C'_s)$,
where $s=1,\ldots,d/2$ and one more component $C_0$. 
We choose numbers $a$ and $b$ as before and consider the points $a(s),b(s),\zeta_j(s)$ (resp., $a'(s),b'(s),\zeta'_j(s)$)
on $C_s$ (resp., $C'_s$), where $s=1,\ldots,d/2$. In addition we are going to define certain points on $C_0$. 
Note that the number $m=2g-2-rd$ is even. Applying Theorem \ref{P1thm} we can find distinct complex numbers
$x_1,\ldots,x_{m/2+1};y_1,\ldots,y_{m/2+1};z_1,\ldots,z_n$ such that the rational function
$$R_0(t)=\frac{\prod_{i=1}^n(t-z_i)^{m_i}}{\prod_{j=1}^{m/2+1}(t-x_j)(t-y_j)}$$
has opposite residues at $x_j$ and $y_j$ for all $j$. Let $x_j(0)$, $y_j(0)$, $z_i(0)$ denote the corresponding
points on $C_0$. To obtain $C$ we make the following identifications:

\noindent (a) $a(s)$ is glued to $b(s)$ for $s=1,\ldots,d/2$;

\noindent (a') $x_j(0)$ is glued to $y_j(0)$ for $j=2,\ldots,m/2+1$;

\noindent (b) for $j=2,\ldots,r$; $s=1,\ldots,d/2$, each $\zeta_j(s)$ is glued to $\zeta'_j(s)$;

\noindent (c) for $s=1,\ldots,d/2-1$, $\zeta_1(s)$ is glued to $\zeta'_1(s+1)$; $\zeta_1(d/2)$ is glued to $y_1(0)$;
$x_1(0)$ is glued to $\zeta'_1(1)$.

In other words, $C$ is a wheel formed by $d/2+1$ curves $C_s\cup C'_s$ (glued pairwise as before) and $C_0$,
where $C_0$ has $m/2$ nodes. We set $D=\sum_{s=1}^{d/2} 0(s)+0'(s)$ as before. Again it is clear that $h^0(D)=1$.
Also we set $p_i=z_i(0)$ for $i=1,\ldots,n$.
It remains to show that $\OO_C(rD+\sum_i m_ip_i)\simeq\om_C$. For this we define a regular $1$-form $\eta$ on $C$ vanishing
on $rD+\sum_i m_ip_i$ by setting $\eta|_{C_s}=R(t)dt$, $\eta|_{C'_s}=-R(t)dt$ and $\eta|_{C_0}=\la\cdot R_0(t)dt$, where
$$\la=\frac{\Res_{t=\zeta_1}R(t)dt}{\Res_{t=x_1}R_0(t)dt}.$$

Finally, consider the case when $d$ is odd. Then we define $C$ by gluing from $d$ copies of $\P^1$ grouped in $(d-1)/2$
pairs $(C_s,C'_s)$, where $s=1,\ldots,(d-1)/2$, and one special component $C_0$. The special points $a(s),b(s),\zeta_j(s)$
and $a'(s),b'(s),\zeta'_j(s)$ on components $C_s$ and $C'_s$ (where $s\ge 1$) are the same as in the previous case.
To construct special points on $C_0$ we observe that the number $q=\sum_i m_i+r=2g-2-r(d-1)$ is even, so
applying Theorem \ref{P1thm} we can find distinct complex numbers
$x_1,\ldots,x_{q/2+1};y_1,\ldots,y_{q/2+1};z_0,z_1,\ldots,z_n$ such that the rational function
$$\wt{R}_0(t)=\frac{(t-z_0)^r\prod_{i=1}^n(t-z_i)^{m_i}}{\prod_{j=1}^{q/2+1}(t-x_j)(t-y_j)}$$
has opposite residues at $x_j$ and $y_j$ for all $j$. As before we consider the corresponding points
$x_j(0)$, $y_j(0)$, $z_i(0)$ on $C_0$ and make the following identifications to get $C$: 

\noindent (a) $a(s)$ is glued to $b(s)$ for $s=1,\ldots,(d-1)/2$;

\noindent (a') $x_j(0)$ is glued to $y_j(0)$ for $j=2,\ldots,q/2+1$;

\noindent (b) for $j=2,\ldots,r$; $s=1,\ldots,d/2$, each $\zeta_j(s)$ is glued to $\zeta'_j(s)$;

\noindent (c) for $s=1,\ldots,(d-1)/2-1$, $\zeta_1(s)$ is glued to $\zeta'_1(s+1)$; $\zeta_1((d-1)/2)$ is glued to $y_1(0)$;
$x_1(0)$ is glued to $\zeta'_1(1)$.

Thus, $C$ is a wheel formed by $C_s\cup C'_s$ and $C_0$ as above, where $C_0$ has $q/2$ nodes. We set 
$D=z_0(0)+\sum_{s=1}^{d/2} 0(s)+0'(s)$ (so now $D$ has one point on each component). Since $q>0$, each component
has at least one node, so $h^0(D)=1$. Also, as before we set $p_i=z_i(0)$ for $i=1,\ldots,n$. The regular $1$-form
$\eta$ on $C$ with the divisor $rD+\sum_i m_ip_i$ is defined in exactly the same way as above with $R_0$ replaced
by $\wt{R}_0$.
\ed

\begin{rem} In the case when $r$ is even there is a simpler choice of a nodal rational curve in the above construction.
For example, if $m_i=0$ for all $i$, we can just take $C$ to be a wheel of $d$ components, where each component
has $r/2$ nodes. If $m=\sum_i m_i>0$ one can take $C$ to be a wheel of $d+1$ components, $d$ of which have $r/2$ nodes
and one special component has $m/2$ nodes.
\end{rem}  

\subsection{Smoothening}

Let us fix some data $(g, r, n, \bm, d)$ as in the introduction, where $d\ge 0$.
Let $\Mgnbar$ be the moduli space of stable $n$-pointed curves of genus $g$,
$\CC^{\reg}\ra\Mgnbar$ be the open part of the universal curve obtained by removing
all nodes. Then the relative symmetric product $(\CC^{\reg})^{(d)}$ parametrizes
collections $(C,D,p_1,\ldots,p_n)$, where $(C,p_1,\ldots,p_n)$ is stable, $D$
is an effective divisor of degree $d$ supported on the smooth part of $C$.
Let $\ZZ\sub(\CC^{\reg})^{(d)}$ be the degeneracy locus corresponding to collections
with $h^0(\om_C(-rD-\sum_i m_ip_i))\neq 0$: we define it using
the morphism of vector bundles (\ref{mainmap}). Let $\ZZ^0\sub\ZZ$ be the open substack
corresponding to collections $(C,D,p_1,\ldots,p_n)$ such that for every irreducible
component $C_i\sub C$ one has 
$$\deg\om_C(-rD-\sum_i m_ip_i)|_{C_i}=0$$
and such that $h^0(D)=1$.
Theorem \ref{nodalthm} states that $\ZZ^0$ is always nonempty.
However, the point constructed in this theorem lives on the boundary (i.e., on the
preimage of the boundary in $\Mgnbar$). To prove part (b) of 
Theorem \ref{mainthm} we have to prove that there exists a point of $\ZZ^0$
with smooth $C$.
The same dimension count as in section \ref{calcsec} shows that
the dimension of $\ZZ^0$ at every point is $\ge 2g-2+n+d$.
Thus, we will be able to find a point of $\ZZ^0$ with smooth $C$ (and therefore
finish the proof of Theorem \ref{mainthm}) once we 
prove the following result.

\begin{thm}\label{boundest} 
Let $\De\sub\ZZ^0$ be the boundary divisor, i.e., the preimage of the boundary in $\Mgnbar$
under the natural morphism $\ZZ^0\ra\Mgnbar$. Then $\dim\De<2g-2+d+n$.
\end{thm}

The idea of the proof is to estimate the dimension over each strata of the boundary using the tangent
space calculations similar to \ref{calcsec}. Namely, let us consider the modified data
$(g,r,n,d,s,\bm)$ consisting of integers $g\ge 1$, $r\ge 2$, $n\ge 0$, $d\ge 0$, $s\ge 1$ 
and a vector $\bm=(m_1,\ldots,m_n)$ of non-negative integers such that
$$m_1+\ldots+m_n+rd=2g-2+s.$$
To these data we can associate the moduli space $\Meffptpt$ parametrizing collections
$(C,D,p_1,\ldots,p_n,q_1,\ldots,q_s)$, 
where $C$ is a (connected) smooth complex projective
curve of genus $g$, $D\sub C$ is an effective divisor of degree $d$,
$p_1,\ldots,p_n,q_1,\ldots,q_s$ are distinct marked points on $C$ such that
$q_1,\ldots,q_s$ do not belong to the support of $D$ and
\begin{equation}
\OO_C(rD+m_1p_1+\ldots+m_np_n)\simeq\om_C(q_1+\ldots+q_s).
\end{equation}
More precisely, we define $\Meffptpt$ from the cartesian diagram
\begin{equation}
\begin{diagram}
\Meffptpt &\rTo& \XX^d \nonumber\\
\dTo &&\dTo^{\wt{\si}_r^{\bm}}\\
\Mgns &\rTo^{c}& \JJ^{2g-2}
\end{diagram}
\end{equation}
where $\CC\ra\Mgns$ is the universal curve, $\XX^d\sub\CC^{(d)}$ is the open
subset parametrizing effective divisors $D$ such that the support of $D$ does not
contain any of the points $q_1,\ldots,q_s$. Also in this diagram $\JJ^{2g-2}$ denotes 
the relative Jacobian of degree $2g-2$ over $\Mgns$ and 
the morphism $\wt{\si}_r^{\bm}:\XX^d\ra\JJ^{2g-2}$
sends $(C,D,p_1,\ldots,p_n,q_1,\ldots,q_s)$ to $(C,\OO_C(rD+\sum_i m_ip_i-\sum_j q_j))$
(as before, the morphism $c$ sends a curve $C$ to $(C,\om_C)$). 
Note that if $s=1$ then from the Residue Theorem we immediately obtain that above
moduli space is empty.

Our calculations in section \ref{calcsec} can be easily modified to prove the following.

\begin{thm}\label{transdimthm} The stack $\Meffptpt$ is smooth of dimension
$2g-3+d+n+s$ at every point.
\end{thm}

\Pf . The above cartesian diagram shows that at every point 
$$\dim\Meffptpt\ge \dim\XX^d-g=2g-3+d+n+s.$$
Thus, it is enough to prove that at every point $(C,D,p_1,\ldots,p_n,q_1,\ldots,q_s)$
of $\Meffptpt$ the map
$$
\begin{diagram}
T\CC^{(d)}\oplus T\Mgns & \rTo^{d\wt{\si}_r^{\bm},dc} & T\JJ^{2g-2}
\end{diagram}
$$
is surjective (where the tangent spaces are taken at induced points of the relevant spaces).
As in section \ref{calcsec} we can identify this map with the map induced on $H^1$
by the natural morphism of sheaves
\begin{equation}\label{mainmodmor}
\TT_C(-D-E)\oplus\TT_C(-E)\ra A_{\om_C,p,q},
\end{equation}
where $E=\sum_i p_i+\sum_j q_j$, $p=(p_1,\ldots,p_n)$, $q=(q_1,\ldots,q_s)$.
As before, the second component of the latter morphism 
is given by the Lie derivative:
$$\TT_C(-E)\ra A_{\om_C,p,q}: v\mapsto L_v$$
where $L_v$ is consideres as an operator from $\om_C$ to itself.
To describe the first component
let us set $L=\OO_C(D)$ and let $s_1\in L$ be the natural section vanishing at $D$.
Then the first component is given by the natural morphism
\begin{equation}\label{1stcomp}
\TT_C(-D-E)\ra A_{L,p,q}: v\mapsto (s\mapsto v(s/s_1)s_1)
\end{equation}
taking into account the isomorphism
$A_{L,p,q}\simeq A_{\om_C,p,q}$ obtained by applying Lemma \ref{rootdeflem} twice.
More precisely, if we denote by $s_0$ the (unique up to a scalar) rational section of $\om_C\otimes L^{-r}$
vanishing along $\sum_i m_i p_i$ and with poles at $\sum_j q_j$, then the isomorphism
$$A_{\om_C,p,q}\ra A_{L,p,q}$$
sends $\del\in A_{\om_C,p,q}$ to the operator 
$$s\mapsto \frac{\del(s^r s_0)}{rs^{r-1}s_0}$$
from $L$ to itself.
As before, morphism (\ref{1stcomp}) fits into the exact sequence
$$0\ra \TT_C(-D-E)\ra A_{L,p,q}\ra L\ra 0$$
where the morphism $A_{L,p,q}\ra L$ sends $\del\in A_{L,p,q}$ to $\del(s_1)\in L$.
Thus, as in the proof of Proposition \ref{cokerprop} the cokernel of the map on $H^1$ induced by (\ref{mainmodmor})
coincides with the cokernel of the map on $H^1$ induced by the composition
$$\TT_C(-E)\ra A_{\om_C,p,q}\simeq A_{L,p,q}\ra L.$$
It is easy to see that this composed map sends $v\in\TT_C(-E)$ to $(d\lan s_1^rs_0,v\ran)/rs_1^{r-1}s_0$.
This map fits into the following morphism of exact sequences:
\begin{equation}
\begin{diagram}[size=5em]
0 \rTo^{\ } &\TT_C(-E) &\rTo^{s_1^r s_0}&\OO_C(-E')&
\rTo &\OO_C(-E')/\OO_C(-rD-E_{\bm}-E')&\rTo 0 \nonumber\\
&\dTo &&\dTo{d} &&\dTo{f} & \\
0\rTo & L &\rTo^{rs_1^{r-1}s_0}&\om_C(E'')&\rTo &\om_C(E'')/\om_C(E''-(r-1)D-E_{\bm})&\rTo 0
\end{diagram}
\end{equation}
where $E'=p_1+\ldots+p_n$, $E''=q_1+\ldots+q_s$ and $E_{\bm}=m_1p_1+\ldots+m_np_n$. 
Since $H^1(\om_C(\sum_j q_j))=0$, our assertion would follow from the surjectivity of the morphism $f$ in
this diagram. But $f$ factors as a composition of the map
\begin{equation}\label{diffmap}
\OO_C(-E')/\OO_C(-rD-E_{\bm}-E')\ra\om_C/\om_C(-(r-1)D-E_{\bm})
\end{equation}
induced by the de Rham differential with the natural isomorphism
$$\om_C/\om_C(-(r-1)D-E_{\bm})\wt{\ra}\om_C(E'')/\om_C(E''-(r-1)D-E_{\bm})$$
(here we use the assumption that the supports of $E''$ and of $D+p_1+\ldots+p_n$ do not intersect).
Now the surjectivity of (\ref{diffmap}) follows Lemma \ref{divdifflem}, hence, $f$ is surjective.
\ed

\noindent
{\it Proof of Theorem \ref{boundest}.}
Let $S\sub\De$ be one of the strata of the boundary divisor obtained by fixing the combinatorial
type of the stable $n$-pointed curve. It is enough to prove that $\dim S<2g-2+n+d$. Let $\Ga$ be
the corresponding graph (vertices of $\Ga$ correspond to irreducible components of a curve in $S$).  
We denote by $v$ and $e$ the number of vertices and egdes in $\Ga$ respectively. Then
for every $(C,D,p_1,\ldots,p_n)\in S$ the genus $\wt{g}$ of the normalization $\wt{C}$ is equal to
$$\wt{g}=g-1+v-e.$$
Note that the normalized curve $\wt{C}$ consists of $v$ connected components and the map
$\wt{C}\ra C$ glues pairwise $2e$ points $q_1,\ldots,q_{2e}$ on $\wt{C}$. 
Moreover, the data $(D,p_1,\ldots,p_n)$ define similar data $(\wt{D},\wt{p}_1,\ldots,\wt{p}_n)$
on $\wt{C}$ such that
$$\OO_{\wt{C}}(r\wt{D}+\sum_i m_i\wt{p}_i)\simeq\om_{\wt{C}}(\sum_{j=1}^{2e} q_j).$$
Restricting these data to connected components of $\wt{C}$ we obtain points of the moduli spaces
considered in Theorem \ref{transdimthm}. Moreover,
it is clear that the point $(C,D,p_1,\ldots,p_n)\in S$ is uniquely recovered from these points.
Therefore, Theorem \ref{transdimthm} implies the following estimate:
$$\dim S\le 2\wt{g}-3v+d+n+2e.$$
Combining this with the above formula for $\wt{g}$ we get
$$\dim S\le 2g-2+d+n-v<2g-2+d+n.$$
\ed

\section{Complements}\label{compsec}

\subsection{Connected components}\label{conncompsec}

Recall that $\HH(m_1,\ldots,m_n)$ denotes the moduli space of pairs
$(C,\om)$ where $\om$ is a nonzero holomorphic $1$-form on a 
smooth compact complex curve $C$
with multiplicities of zeroes $(m_1,\ldots,m_n)$ (see \cite{KZ}).
Let also $\HH^{num}(m_1,\ldots,m_n)$ be the moduli space of holomorphic
$1$-forms on curves with {\it numbered} zeroes such that the first zero
has multiplicity $m_1$ etc. There is a natural finite morphism
$\HH^{num}(m_1,\ldots,m_n)\ra \HH(m_1,\ldots,m_n)$ inducing a bijection
between the sets of connected components
(see Remark 1 on p.632 of \cite{KZ}).

\noindent
{\it Proof of Theorem \ref{conncompthm}.} Let us consider the stack
$$\HH^{num}(r^d,m_1,\ldots,m_n):=\HH^{num}(r,\ldots,r,m_1,\ldots,m_n)$$
where $r$ is repeated $d$ times in the right hand side.
We have a natural action of the symmetric group $S_d$ on this stack and
a natural finite map to $\HH(r^d,m_1,\ldots,m_n)$ that 
factors as the following composition
$$\HH^{num}(r^d,m_1,\ldots,m_n)\ra \HH^{num}(r^d,m_1,\ldots,m_n)/S_d\ra 
\HH(r^d,m_1,\ldots,m_n).$$
Since the composed map induces a bijection on the sets of
connected components 
we deduce that $\HH^{num}(r^d,m_1,\ldots,m_n)/S_d$ has the same set
of connected components as $\HH(r^d,m_1,\ldots,m_n)$.

Now let us consider the natural map
$$\HH^{num}(r^d,m_1,\ldots,m_n)/S_d\ra\Meffpt$$
sending
the differential $\om$ with divisor of zeroes 
$rq_1+\ldots+rq_d+m_1p_1+\ldots+mp_n$
to $(D,p_1,\ldots,p_n)$, where $D=q_1+\ldots+q_d$.
The image of this map is
the open set $U\sub\Meffpt$ consisting of $(C,D,p_1,\ldots,p_n)$
such that the divisor $D$ is simple (i.e., $D=q_1+\ldots+q_d$ for distinct
points $q_1,\ldots,q_d$)
 and the points $p_i$ are disjoint from
the support of $D$. Since $\HH^{num}(r^d,m_1,\ldots,m_n)/S_d$
is a $\C^*$-torsor over $U$, it remains to prove that the complement
to $U$ in $\Meffpt$ has codimension $\ge 1$.
A point $(C,D,p_1,\ldots,p_n)\in\Meffpt$ belongs to this complement if
either the divisor $D$ is not simple or some of the points $p_i$ belong
to the support of $D$. Thus, $\Meffpt\setminus U$
has a natural stratification corresponding to fixing multiplicities
of $D$ and the pattern according to which some of the points 
$p_i$ collide with points in the support of $D$. 
Over each of these strata there is a natural $\C^*$-torsor that
can be identified with a quotient by a finite group action 
of one of the spaces $\HH(k_1,\ldots,k_N)$ where $N<d+n$. Since
the dimension of $\HH(k_1,\ldots,k_N)$ is equal to $2g-1+N$ this
implies that $\dim(\Meffpt\setminus U)<2g-2+d+n$.  
\ed

Combining Theorem \ref{conncompthm} with Theorems 1 and 2 of \cite{KZ}
we obtain the following information about connected components of our
moduli spaces.

\begin{cor}\label{conncor1} 
Assume that $g\ge 4$.
The moduli space $\Meffpt$ has 
three connected components in the following cases:

\noindent
(i) $r\ge 2$, one marked point, $m_1=2g-2$;

\noindent
(ii) $r\ge 2$, two marked points, $m_1=m_2=g-1$, $g$ is odd;

\noindent
(iii) $r=2g-2$, no marked points;

\noindent
(iv) $r=g-1$, one marked point, $m_1=g-1$, $g$ is odd;

\noindent
(v) $r=g-1$, no marked points, $g$ is odd.

If either all the numbers $r, m_1,\ldots,m_n$ are even or 
all $m_1,\ldots,m_n$ are even and $\sum m_i=2g-2$ (so that $d=0$), 
then $\Meffpt$ has two connected components. In addition $\Meffpt$ has
two connected components in the following cases:

\noindent
(a) $r\ge 2$, two marked points, $m_1=m_2=g-1$, $g$ is even;

\noindent
(b) $r=g-1$, one marked point, $m_1=g-1$, $g$ is even;

\noindent
(c) $r=g-1$, no marked points, $g$ is even.

All the other moduli spaces $\Meffpt$ for $g\ge 4$ are connected.
\end{cor}

\begin{cor}\label{conncor2}
All the spaces $\Mtwoeffpt$ are connected.

The space $\Mteffpt$ has two connected components in the following cases:

\noindent
(i) $r\ge 2$, one marked point, $m_1=4$;

\noindent
(ii) $r\ge 2$, two marked points, $m_1=m_2=2$;

\noindent
(iii) $r=4$, no marked points;

\noindent
(iv) $r=2$, one marked point, $m_1=2$;

\noindent
(v) $r=2$, no marked points.

All the other spaces $\Mteffpt$ are connected.
\end{cor}

\begin{ex}
Let us consider our moduli spaces in the case $r=2$, no marked points,
$g\ge 3$. The corresponding spaces
$\Mtwoeff$ parametrize pairs $(C,D)$ where $D$ is an effective divisor
of degree $g-1$ on $C$ such that $2D$ is a canonical divisor. 
Theorem \ref{conncompthm} implies that
$\Mtwoeff$ has two connected components, that are distinguished by
the parity of $h^0(D)$. The image of the even component under the forgetting
map to $\MM_g$
is the divisor $\DD\sub\MM_g$ consisting of curves posessing a singular
even theta-characteristic, i.e., a square root $L$ of the canonical bundle
such that $h^0(L)$ is even and $h^0(L)>0$. Thus, we obtain a new proof
of connectedness of $\DD$ that was established by Teixidor i Bigas in 
\cite{T2}. 
On the other hand, Theorem 2.16 of \cite{T1} asserts that the curve 
corresponding to a generic point in $\DD$ has only one even 
theta-characteristic $L$ with $h^0(L)\ge 2$. 
Thus, the even component of $\Mtwoeff$ can be viewed as a resolution of
singularities of $\DD$.
\end{ex}

\subsection{Witten's top Chern class}

We can calculate the difference $c^{\frac{1}{r}}-c_{-\chi}(R\pi_*\LL)$ in terms
of the locus $S^0$ in the following examples. We adopt the conventions
of \cite{V} regarding the intersection theory on Deligne-Mumford stacks.

\noindent
1. $r=2$, no points, $g\ge 1$. The class $c^{\frac{1}{2}}$ has degree zero.
More precisely, the moduli space $\MM=\MM_g^{\frac{1}{2}}$ has two components
$\MM^+$ and $\MM^-$ corresponding to even and odd theta-characterstics and
we have 
$$c^{\frac{1}{2}}=[\MM^+]-[\MM^-]$$
(see \cite{PV}, 5.3). On the other hand, clearly
$c_0(R\pi_*\LL)=[\MM^+]+[\MM^-]$. 
The only component of $S^0$ of maximal dimension is $S^{0,-}=\MM^-$.
Hence, in this case we have
$$c^{\frac{1}{2}}-c_0(R\pi_*\LL)=-2[S^{0,-}].$$

\noindent
2. $g=1$, no points, $r\ge 2$. The class $c^{\frac{1}{r}}$ still has
codimension zero. Again using the known formula for $c^{\frac{1}{r}}$ in this case (see
\cite{PV}, 5.3) we find
$$c^{\frac{1}{r}}-c_0(R\pi_*\LL)=-r[S^0].$$

\noindent
3. $g=2$, one marked point with $m_1=2$, $r\ge 2$. 
In this case the restriction of the projection $\Mtospin\ra\Mto$ to
the locus $S^0$ gives an isomorphism of $S^0$
with the divisor of pairs $(C,p)$ such that $p$ is a
Weierstrass point on $C$, i.e., $\om_C\simeq\OO_C(2p)$. The irreducibility of
this divisor is well known (and also follows from Corollary \ref{conncor2}).
The class $c^{\frac{1}{r}}$ in this case has degree $1$ and we want to compare it with
the standard divisor class $\mu_1:=c_1(R\pi_*\LL)$, i.e., 
compute the coefficient $m$ such that
$$c^{\frac{1}{r}}-\mu_1=m[S^0],$$
where $\MM=\Mtospin$.
Let us consider the closure $\ov{S}^0\sub\Mtospinbar$ of $S^0$.
Then it is easy to see that the locus in $\Mtospinbar$ where the spin structure has a non-zero
global section is the union of $\ov{S}^0$ with an irreducible boundary divisor $\de_0\sub\Mtospinbar$
containing curves glued from two elliptic curves.
Therefore, we have
\begin{equation}\label{meq}
c^{\frac{1}{r}}-\mu_1=m[\ov{S}^0]+m'[\de_0]
\end{equation}
for some $m'$. 
Now the idea is to use the Ramond vanishing axiom for $c^{\frac{1}{r}}$ stating that
the restriction of $c^{\frac{1}{r}}$ to the Ramond component of the boundary in $\Mtospinbar$ vanishes 
(see \cite{P}).
Namely, we are going to construct a family of $r$-spin curves over a $1$-dimensional stack $B$ that corresponds
to a morphism from $B$ to the Ramond boundary component and then calculate the pull-back to $B$ of
the divisor classes appearing in the above equality.

Let us fix an elliptic curve $E$ without complex multiplication and a pair of distinct points $p,q\in E$ such that
$2p$ is not rationally equivalent to $2q$. 
We denote by $C$ the nodal curve of arithmetic genus $2$
obtained from $E$ by gluing $p$ with $q$. The family over $E\setminus\{p,q\}$ that associates
to every point $x$ the $1$-pointed curve $(C,x)$ has a natural completion to a family of stable
curves over $E$. Namely, at points $p$ and $q$ the curve $C$ 
gets replaced by the curve that is obtained by gluing $E$ with $\P^1$ along $2$ points that correspond
to $p$ and $q$ on $E$ (say, $p$ is glued to $0\in\P^1$ and $q$ is glued to $\infty\in\P^1$)
and the marked point belongs to $\P^1$. Here is a more precise definition of this family.
Let $\wt{\CC}\ra E\times E$ be the blow-up of $E\times E$
at points $(p,p)$ and $(q,q)$. Let $\De'\sub\wt{\CC}$, $D_p$ and $D_q$ be the proper preimages of the
diagonal $\De$ and of the divisors $E\times\{p\}$ and $E\times\{q\}$ in $E\times E$. 
Let $p_1:\wt{\CC}\ra E$ be the first projection. Note that all three divisors $\De'$, $D_p$ and $D_q$
are sections of $p_1$, hence $(p_1:\wt{\CC}\ra E, \De',D_p,D_q)$
can be viewed as a family of $3$-pointed curves. Now our family $\pi:\CC\ra E$ is obtained
from $\wt{\CC}$ by gluing the $E$-points $D_p$ and $D_q$ (with $\De'$ viewed as a marked point).
Let $E\ra\Mtobar$ be the morphism corresponding to this family. By definition it factors as a
composition:
$$E\ra \Motbar\ra\de_{\irr}\sub\Mtobar,$$ 
where $\de_{\irr}$ is the component of the boundary divisor
in $\Mtobar$ containing irreducible curves of geometric genus $1$ with one node, $\Motbar\ra\de_{\irr}$ is the standard
gluing morphism. Note that the morphism $E\ra\Motbar$ corresponds to the family $\wt{\CC}$ over $E$.
Let us define the $1$-dimensional stack $B$ so that the following diagram is cartesian:
\begin{equation}
\begin{diagram}
B &\rTo{}& \de_{\irr,R} \nonumber\\
\dTo{f} &&\dTo\\
E &\rTo& \de_{\irr}
\end{diagram}
\end{equation}
where $\de_{\irr,R}\sub\Mtospinbar$ is the irreducible component of the boundary divisor in $\Mtospinbar$ containing
irreducible curves of geometric genus $1$ equipped with $r$-spin structures of Ramond type at the node.
The covering $B\ra E$ has degree $r^2$ and is ramified over $p$ and $q$.
Note that we also have the following commutative (non-cartesian) diagram:
\begin{equation}
\begin{diagram}
B &\rTo{}& \Motspinbar
\nonumber\\
\dTo{f} &&\dTo{p}\\
E &\rTo& \Motbar
\end{diagram}
\end{equation}
where the upper horizontal arrow is the composition
$$B\ra\de_{\irr,R}\times_{\de_{\irr}}\Motbar
\ra\Motspinbar.$$
Here we use the standard gluing morphism for the Ramond boundary (see \cite{JKV}).
It is easy to see that the pull-backs to $\de_{\irr,R}\times_{\de_{\irr}}\Motbar$
of the class $\mu_1$ on the moduli space $\Mtospinbar$ and the class $\mu_1$ on $\Motspinbar$
coincide (this is always true for the Ramond gluing). Now on the moduli space 
$\Motspinbar$ we can use the formula for the class $\mu_1$ given in Proposition 2.4 of \cite{JKV}:
\begin{equation}\label{mu1eq}
2r^2\mu_1=(2r^2-12r+12)p^*\la_1-r\eps+(r-1)\de+\sum_i m_i(r-2-m_i)p^*\psi_i
\end{equation}
where $\la_1$ and $\psi_i$ ($i=1,2,3$) are standard classes on $\Motbar$, $\eps$ and $\de$ are certain
combinations of the boundary divisors on $\Motspinbar$. Note that the image of $B$ intersects only with
one boundary divisor $\wt{\a}$ in $\Motspinbar$. The generic point of $\wt{\a}$ corresponds to
a reducible $3$-pointed curve $C_0\cup C_1$ with one node, where $C_0$ has genus $0$, $C_1$ has genus $1$;  
the first two marked points (with markings $2$ and $r-1$) belong to $C_0$, the third marked
point (with the marking $r-1$) belongs to $C_1$. The type of the $r$-spin structure at the node in this
case is determined uniquely: one has $(m^+,m^-)=(r-3,1)$. Hence, by Proposition 2.1 of \cite{JKV} we have
$$\wt{\a}=\frac{gcd(2,r)}{r}p^*\a$$
where $\a$ is the boundary divisor on $\Motbar$ containing reducible $3$-pointed curves as above (with the first two marked points 
on the rational component); $gcd(m,n)$ denotes the greatest common divisor of $m$ and $n$.
Also, using the definition in Proposition 2.4 of \cite{JKV} we find
$$\eps=\frac{2(r-2)}{gcd(2,r)}\wt{\a}+\ldots,$$
$$\de=\frac{r}{gcd(2,r)}\wt{\a}+\ldots$$
where the dots denote combinations of other boundary components.
It follows that 
$$\eps\cdot[B]=\frac{2(r-2)}{r}p^*\a\cdot[B]=\frac{2(r-2)}{r}f^*(\a\cdot[E]);$$
$$\de\cdot[B]=p^*\a\cdot[B]=f^*(\a\cdot[E]),$$
where of course $\a\cdot[E]=p+q$.
It is also easy to calculate that
$$\la_1\cdot[E]=0,\ \psi_1\cdot[E]=p+q,\ \psi_2\cdot[E]=p,\ \psi_3\cdot[E]=q.$$
Hence, using the above formula for $\mu_1$ we get
$$2r^2\mu_1\cdot[B]=[-2(r-2)+r-1+2(r-4)-(r-1)]f^*(p+q)=-4f^*(p+q),$$
therefore,
$$\mu_1\cdot[B]=-\frac{2}{r^2}f^*(p+q).$$

It remains to compute $[\ov{S}^0]\cdot[B]$.
It is easy to see that the support of this divisor does not intersect $f^{-1}(\{p,q\})$,
hence it consists of $4$ points in $B$ that project to $4$ points $x\in E$ such that
$2x$ is rationally equivalent to $p+q$. We claim that all these points have multiplicity $1$
in $[\ov{S}^0]\cdot[B]$. To prove this let us consider the natural projection
$\phi:\ov{S}^0\ra Z$, where $Z\sub\Mtobar$ is the locus of curves $(C,x)$ with 
$h^0(\om_C(-2x))\neq 0$ defined as the degeneracy locus of the map $H^0(\om_C)\ra H^0(\om_C|_{2x})$.
Let $U\sub\Mtobar$ (resp., $U^{\frac{1}{r}}\sub\Mtospinbar$) be the complement in $\Mtobar$ (resp.,
$\Mtospinbar$) to the union of all boundary components different from $\de_{\irr}$ (resp.,
$\de_{irr,R}$). Then it is easy to see that the restriction of $\phi$ induces an isomorphism
$$\ov{S}^0\cap U^{\frac{1}{r}}\ra Z\cap U.$$
The point is that over $U$ the condition $h^0(\om_C(-2x))\neq 0$ implies $\om_C(-2x)\simeq\OO_C$,
while over $U^{\frac{1}{r}}$ the $r$-spin structure $L$ is locally free and the condition that
$h^0(L)\neq 0$ implies $L\simeq\OO_C$. Therefore, to prove our claim it suffices to
check that the points $x\in E$ such that
$2x\sim p+q$, have multiplicity one in $[Z]\cdot[E]$. Note that the support of $[Z]\cdot[E]$ consists of
the points $p$, $q$ and $4$ points $x$ with $2x\sim p+q$. On the other hand, from the definition
of $Z$ as the degeneracy locus we have $[Z]\sim 3\psi_1-\la$. Hence, the divisor $[Z]\cdot[E]$ is rationally
equivalent to $3(p+q)$. Therefore, $\deg[Z]\cdot[E]=6$ which implies that all points in the support of
$[Z]\cdot[E]$ have multiplicity $1$.

The above argument implies that
$$f_*([\ov{S}^0]\cdot[B])=\frac{2}{r}(p+q).$$
On the other hand, we have $f_*(\mu_1\cdot[B])=-2(p+q)$. Therefore, restricting the equality (\ref{meq})
to $B$ and then pushing to $E$ (taking into account the vanishing of $c^{\frac{1}{r}}\cdot[B]$) we find
that $m=r$, i.e.,
\begin{equation}\label{divclasseq}
c^{\frac{1}{r}}-\mu_1=r[S^0].
\end{equation}

\noindent
4. $g=2$, two marked points, $m_1=m_2=1$, $r\ge 2$. In this case one can show 
that the equation (\ref{divclasseq}) still holds. Namely, we can fix an elliptic curve $E$
with $3$ distinct points $p,q,x_0$ such that $2x_0\not\sim p+q$ and consider the nodal curve $C$ of 
arithmetic genus $2$ obtained by identifying $p$ and $q$. Then we can consider the family of $2$-pointed
curves with the base $E$ whose generic member is $(C,x_0,x)$, where $x\in E$. As above we can construct
a family of $r$-spin curves with the stacky base $B$ given by some covering of $E$, such that
the $r$-spin structures in this family are of Ramond type near one node. The rest of the argument
is very similar to the previous case so we skip it.

\noindent
5. $g=3$, one marked point, $m_1=1$, $r=3$. Again we claim that the equation (\ref{divclasseq}) holds in this case.
One can check this by considering the family parametrized by a smooth curve $C$ of genus $2$.
Namely, we fix two generic points $p,q\in C$ and consider the nodal curve $C_0$ obtained by gluing together
$p$ and $q$. Varying a point $x$ on $C$ we get a family of nodal curves with generic member $(C_0,x)$.
Then as above we consider the covering $f:B\ra C$ corresponding to adding the $r$-spin structures.
The computation of $f_*(\mu_1\cdot[B])$ in this case is similar to the previous cases. The answer turns out to be
$$f_*(\mu_1\cdot[B])=-9(p+q).$$
The class $f_*([\ov{S}^0]\cdot[B])$ is equal to $[D]/3$, where $D$ is the following divisor on $C$:
$$D=\sum_{y:\om_C(p+q-3y)\simeq\OO_C(x)}x.$$
It can be computed as follows. Note that the divisor 
$$D'=\sum_{y:h^0(\om_C(p+q-3y))\neq 0} y$$ 
on $C$ can be presented
as the degeneracy locus of the morphism $H^0(\om_C(p+q))\ra \om_C|_{3y}$. Computing the determinant, we
conclude that $\OO_C(D')\simeq\om_C^6(3p+3q)$. In particular, $\deg(D')=18$. Therefore,
$$\OO_C(D)\simeq\om_C^{18}(18p+18q-3D')\simeq\OO_C(9p+9q).$$
Comparing this with the formula fo $f_*(\mu_1\cdot[B])$ we get the result.

\noindent
6. $g\ge 3$, one marked point, $m_1=2$, $r=2$. We claim that the equation (\ref{divclasseq}) still holds.
Note that the moduli space $\MM=\MM_{g,1}^{\frac{1}{2},2}$ in this case has two connected components
$\MM^+$ and $\MM^-$
depending on the parity of $h^0(L(p))$, where $L$ is a square root of $\om_C(-2p)$ ($p$ is a marked point). 
Let us consider the situation
over these two components separately. Over the complement to a closed subset of codimension $\ge 2$ in
the component $\MM^-$ we have $h^0(L(p))=1$ and the locus $S^0$
coincides with the zero locus of the map $H^0(L(p))\ra L(p)|_p$. 
Using the isomorphism
$$H^0(L(p))^{\otimes 2}\simeq\det R\Ga(L(p))\simeq\det R\Ga(L)\otimes L(p)|_p$$
one can easily deduce that
$$[S^0]=\frac{1}{2}(\wt{\psi}-\mu_1)$$
where $\wt{\psi}=\psi_1/2$ is the class of the line bundle with the fiber $L(p)|_p$ over $(C,p,L)\in\Mtospintwo$.
But the descent axiom (Proposition 5.1 of \cite{PV}
\footnote{The class considered in \cite{PV} differs from our $c^{\frac{1}{r}}$ by the sign $(-1)^{\chi(L)}$.}) 
implies that over $\MM^-$ one has
$c^{\frac{1}{2}}=\wt{\psi}$, so the equation (\ref{divclasseq}) follows. Now consider the situation over the even
component $\MM^+\sub\MM$. By the descent axiom in this case we have $c^{\frac{1}{2}}=-\wt{\psi}$.
On the complement to a closed subset of codimension $\ge 2$ in $\MM^+$ we have $h^0(L(p))\le 2$ and $h^0(L)\le 1$.
Hence, over this complement $S^0$ parametrizes data $(C,p,L)$ such that $L(p)$ is an even theta-characteristic
and $h^0(L(p))=2$. Let $\phi:\MM^+\ra\MM_{g,1}$ be the natural projection. Note that
$\phi$ has degree $2^{g-2}(2^g+1)$ (since there are $2^{g-1}(2^g+1)$ even theta-characteristics on a curve of genus $g$).
Then we have 
$$\phi_*[S^0]=\frac{1}{2}[\DD]=2^{g-4}(2^g+1)\la_1,$$ 
where $\DD\sub\MM_{g,1}$ is the divisor of curves having an even theta-characteristic with at least two independent sections
(the formula for $[\DD]$ can be found in \cite{T2}).
On the other hand, using the formula (\ref{mu1eq}) we find 
$$\phi_*(c^{\frac{1}{2}}-\mu_1)=\phi_*(\frac{1}{2}\la_1)=2^{g-3}(2^g+1)\la_1=2\phi_*[S^0].$$

The above examples suggest that the equation (\ref{divclasseq}) holds whenever $\deg c^{\frac{1}{r}}=-\chi=1$.
This can happen only if either $r=2$ or $g\le 4$. More precisely, there are three
remaining cases to check: 1)$r=2$, $g\ge 3$, $m_1=m_2=1$; 2) $g=3$, $r=4$; and 3)$g=4$, $r=3$. More important problem is to find a 
conceptual explanation of the equation (\ref{divclasseq}) and its 
generalization to the case $-\chi\ge 2$.

\end{document}